\documentclass[11pt,a4paper]{article}
\usepackage{amscd}
\usepackage[latin1]{inputenc}
\usepackage{amsmath}
\usepackage{amsfonts}
\usepackage{amssymb}
\usepackage{color}
\usepackage{float}
\usepackage{amsthm}

\usepackage{listings}
\newtheorem{theorem}{Theorem}

\newtheorem{lemma}[theorem]{Lemma}
\newtheorem{proposition}[theorem]{Proposition}
\newtheorem{remark}[theorem]{Remark}

\def\Qed{\hfill\raisebox{.6ex}{\framebox[2.5mm]{}}\\[.15in]}

\def\m{\mathbb}

\begin{document}
\date{}
\title{New canonical triple covers of surfaces}
\author{Carlos Rito}
\maketitle

\begin{abstract}
We construct a surface of general type with canonical map of degree $12$ which factors as
a triple cover and a bidouble cover of $\m P^2$.
We also show the existence of a smooth surface with $q=0,$ $\chi=13$ and $K^2=9\chi$ such that
its canonical map is either of degree $3$ onto a surface of general type or
of degree $9$ onto a rational surface.

\end{abstract}

\section{Introduction}

Let $S$ be a smooth minimal surface of general type.
Denote by ${\phi:S\dashrightarrow\m P^{p_g-1}}$ the canonical map and let $d:=\deg(\phi).$
The following Beauville's result is well-known.
\begin{theorem}[\cite{Be}]
If the canonical image $\Sigma:=\phi(S)$ is a surface, then either:
\begin{description}
\item{{\rm (i)}} $p_g(\Sigma)=0,$ or
\item{{\rm (ii)}} $\Sigma$ is a canonical surface $($in particular $p_g(\Sigma)=p_g(S))$.
\end{description}
Moreover, in case {\rm (i)} $d\leq 36$ and in case {\rm (ii)} $d\leq 9.$
\end{theorem}

Beauville has also constructed families of examples with $\chi(\mathcal O_S)$ arbitrarily
large for $d=2, 4, 6, 8$ and $p_g(\Sigma)=0.$ Although this is a classical problem, for $d>8$ the
number of known examples drops drastically: only Tan's example \cite[\S 5]{Ta2} with
$d=K\sp 2=9,$ $\chi=4$ and Persson's example \cite{Pe} with $d=K\sp 2=16,$ $\chi=4$ are known.
More recently, Du and Gao \cite{DuGa} claim that if the canonical map is an abelian cover
of $\m P^2,$ then these are the only possibilities for $d>8.$

In this note we construct an example with $d=K\sp 2=12$ which factors as
a triple cover and a bidouble cover of $\m P^2.$

Known examples for case {\rm (ii)} with $d=3$ date to 1991/2:
Pardini's example \cite{Pa2} with $K^2=27,$ $\chi=6$ and
Tan's examples \cite{Ta2} with $K^2\leq 6\chi,$ $5\leq\chi\leq 9$.
Nowadays this case $d=3$ is still mysterious. On the one hand no one has given a bound for $\chi,$ on the other hand there are no examples for other values of the invariants.

More generally, for the case where the canonical map factors through a triple cover
of a surface of general type and $d=6,$ we have only the family given in \cite[Example 3.4]{CiPaTo}
with invariants on the Noether's line $K\sp 2=2p_g-4.$

Here we show the existence of a smooth regular surface $S$ with $\chi=13$ and $K^2=9\chi$
such that its canonical map $\phi$ factors through a triple cover of a surface of general type.
This is the first example on the border line $K^2=9\chi$ (recall that, for a surface of
general type, one always has $2p_g-4\leq K^2\leq 9\chi$).

This surface is an unramified cover of a fake projective plane.
We show that if $\deg(\phi)\ne 3,$ then $\phi$ is of degree $9$ onto a rational surface.
If this is the case, then one might expect to be able to recover the construction of the
rigid surface $S$ has a covering of $\m P^2,$ which would be interesting. Since it seems
very difficult to provide a geometric construction of a fake projective plane,
we conjecture that $d=3.$

\bigskip
\noindent{\bf Notation}

We work over the complex numbers. All varieties are assumed to be projective algebraic.
A $(-n)$-curve on a surface is a curve isomorphic to $\m P^1$ with self-intersection $-n.$
Linear equivalence of divisors is denoted by $\equiv.$
The rest of the notation is standard in Algebraic Geometry.\\

\bigskip
\noindent{\bf Acknowledgements}

The author wishes to thank Margarida Mendes Lopes, Sai-Kee Yeung, Gopal Prasad,
Donald Cartwright, Tim Steger and specially Amir Dzambic and Rita Pardini for useful correspondence.

The author is a member of the Center for Mathematics of the University of Porto and is a collaborator of the Center for Mathematical Analysis, Geometry and Dynamical Systems (IST/UTL).
This research was partially supported by the European Regional Development Fund through the programme COMPETE and by the Portuguese Government through the FCT--Funda\c c\~ao para a Ci\^encia e a Tecnologia under the projects PEst--C/MAT/UI0144/2013 and PTDC/MAT-GEO/0675/2012.

\section{Basics on Galois triple covers}
Our references for triple covers are \cite{Ta1}, \cite{Ta2} or \cite{Mi}.

Let $X$ be a smooth surface. A Galois triple cover $\pi:Y\rightarrow X$ is determined by
divisors $L,$ $M,$ $B$, $C$ on $X$ such that $B\in |2L-M|$ and $C\in |2M-L|$.
The branch locus of $\pi$ is $B+C$ and $3L\equiv 2B+C,$ $3M\equiv B+2C$. The surface $Y$
is normal iff $B+C$ is reduced. The singularities of $Y$ lie over the singularities of $B+C$.

If $B+C$ is smooth, we have
\begin{equation}\label{eq1}
\chi(\mathcal O_Y)=3\chi(\mathcal O_X)+\frac{1}{2}\left(L^2+K_XL\right)+\frac{1}{2}\left(M^2+K_XM\right),
\end{equation}
\begin{equation}\label{eq2}
K_Y^2=3K_X^2+4\left(L^2+K_XL\right)+4\left(M^2+K_XM\right)-4LM,
\end{equation}
\begin{equation}\label{eq3}
q(Y)=q(X)+h^1(X,\mathcal O_X(K_X+L))+h^1(X,\mathcal O_X(K_X+M)),
\end{equation}
\begin{equation}\label{eq4}
p_g(Y)=p_g(X)+h^0(X,\mathcal O_X(K_X+L))+h^0(X,\mathcal O_X(K_X+M)).
\end{equation}

Now suppose that $\sigma:X\rightarrow X'$ is the minimal resolution of a normal surface $X'$ with a set $s=\{s_1,\ldots,s_n\}$ of ordinary cusps (singularities of type $A_2$).
If the $(-2)$-curves $A_i,$ $A_i'$ satisfying $\sigma^{-1}(s_i)=A_i+A_i'$ can be
labelled such that
$$\sum_1^n (2A_i+A_i')\equiv 3J,$$
for some divisor $J,$ then we say that $s$ is a {\it $3$-divisible} set of cusps.

\begin{proposition}\label{prop1}
Let $X'$ be a minimal surface of general type containing a $3$-divisible set as above as only singularities. Let $\phi:Y\rightarrow X$ be a Galois triple cover with branch locus $\sum_1^n (A_i+A_i')$.

If $n=3\chi(\mathcal O_X),$ then $\chi(\mathcal O_Y)=\chi(\mathcal O_X)$ and $K_{Y'}^2=3K_X^2,$ where $Y'$ is the minimal model of $Y.$
\end{proposition}
\noindent{\bf Proof:}

Let $\widetilde X\rightarrow X$ be the blow-up at the singular points of $\sum_1^n(A_i+A_i').$ Denote by $\widehat{A_i},\widehat{A_i'}$ the $(-3)$-curves which are the strict transforms of $A_i,A_i',$ $i=1,\ldots,n.$
The surface $Y'$ is the minimal model of the Galois triple cover of $\widetilde X$ with branch locus $\sum_1^n(\widehat{A_i}+\widehat{A_i'}).$
The result follows from (\ref{eq1}) and (\ref{eq2}) (notice that $K_{\widetilde X}^2=K_X^2-n$ and $K_{Y'}^2=3K_{\widetilde X}^2+3n$).

\begin{remark}
Note that the cusps induce smooth points on the covering surface, i.e. the pullback of the divisor $\sum_1^n (A_i+A_i')$ is contracted to smooth points of $Y'.$
\end{remark}

\section{A surface with canonical map of degree $12$}

The following result has been shown by Tan \cite[Thm 6.2.1]{Ta3}, using codes.
Here we give an alternative proof.

\begin{lemma}\label{lemma}
Let $X$ be a double cover of $\m P^2$ ramified over a quartic curve with $3$ cusps.
Then the $3$ cusps of $X$ are $3$-divisible.
\end{lemma}

\noindent{\bf Proof:}
Let $B\subset\m P^2$ be a quartic curve with $3$ cusps at points $p_1,$ $p_2,$ $p_3$
(it is well known that such a curve exists; it is unique up to projective equivalence).
Consider the canonical resolution $X'\rightarrow X.$ The strict transform of the lines
through $p_1p_2,$ $p_2p_3$ and $p_1p_3$ is an union of disjoint $(-1)$-curves
$E_i,E_i'\subset X',$ $i=1,2,3.$ These curves and the $(-2)$-curves $A_i,A_i',$ $i=1,2,3,$ which
contract to the cusps of $X$ can be labelled such that the intersection matrix of the
curves $A_1,A_1',A_2,A_2',A_3,A_3',E_1,E_2$ is
\begin{verbatim}
[-2  1  0  0  0  0  1  0]
[ 1 -2  0  0  0  0  0  0]
[ 0  0 -2  1  0  0  1  0]
[ 0  0  1 -2  0  0  0  1]
[ 0  0  0  0 -2  1  0  0]
[ 0  0  0  0  1 -2  0  1]
[ 1  0  1  0  0  0 -1  0]
[ 0  0  0  1  0  1  0 -1]
\end{verbatim}
This matrix has determinant zero. Since
$$b_2(X')=12\chi(\mathcal O_{X'})-K_{X'}^2+4q(X')-2=8,$$ these $8$ curves are dependent in
${\rm Num} (X'),$ and this relation has to be expressed in the nullspace of the matrix.
Using computer algebra we get that this nullspace has basis
\begin{verbatim}
( 2  1  1 -1 -1 -2  3 -3).
\end{verbatim}
Thus $2A_1+A_1'+A_2-A_2'-A_3-2A_3'+3E_1-3E_2=0$ in ${\rm NS} (X')$ (notice that
$X'$ has no non-trivial torsion).
One has ${\rm NS} (X')={\rm Pic} (X')$ for regular surfaces $X'$ (Castelnuovo),
hence there exists a divisor $L$ such that
$$2A_1+A_1'+A_2+2A_2'+2A_3+A_3'\equiv 3L.$$\Qed

Let $Q_1, Q_2\subset\m P^2$ be quartic curves with $3$ cusps each such that $Q_1+Q_2$
has $6$ cusps and $16$ nodes. Let $V$ be the bidouble cover of $\m P^2$ defined by the divisors
$Q_1,$ $Q_2,$ $Q_3:=0$ (for information on bidouble covers see e.g. \cite{Ca} or \cite{Pa1}).
Consider the divisors $J_1,$ $J_2,$ $J_3$ such that $2J_1\equiv Q_2+Q_3,$ $2J_2\equiv Q_1+Q_3,$ $2J_3\equiv Q_1+Q_2.$
We have $$p_g(V)=p_g\left(\m P^2\right)+\sum_1^3 h^0\left(\m P^2,K_{\m P^2}+J_i\right)=3,$$
$$\chi(\mathcal O_V)=4\chi\left(\mathcal O_{\m P^2}\right)+\frac{1}{2}\sum_1^3 J_i\left(K_{\m P^2}+J_i\right)=4$$
and $V$ has $12$ ordinary cusps and no other singularities.

Denote by $W_1,$ $W_2,$ $W_3$ the double covers of $\m P^2$ with branch curves $Q_2,$ $Q_1,$ $Q_1+Q_2,$ respectively (the intermediate surfaces of the bidouble cover).
The canonical map of $V$ factors through maps
$V\rightarrow W_i,$ $i=1,2,3,$ hence it is of degree $4.$
We get from Lemma \ref{lemma} that the $3$ cusps of $W_1$ and the $3$ cusps of $W_2$
are $3$-divisible, therefore the $12$ cusps of $V$ are also $3$-divisible.

Now let $S\to V$ be the Galois triple cover ramified over the $12$ cusps.
We {\bf claim} that $p_g(S)=p_g(V)=3.$
This implies that the canonical map of $S$ factors through the triple cover,
thus it is of degree $12.$ From Proposition \ref{prop1}, $q(S)=0$ and $K_S^2=12.$

So it remains to prove the claim. Let $\widetilde V$ be the smooth minimal resolution of $V.$
The cusps of $V$ correspond to configurations of $(-2)$-curves $A_i+A_i'\subset\widetilde V,$ $i=1,\ldots,12.$
These can be labelled such that there exist divisors $L,$ $M$ satisfying
$2B+C\equiv 3L,$ $B+2C\equiv 3M,$ where $B:=\sum A_i$ and $C:=\sum A_i'.$

Below we use the notation $D\geq 0$ for
$h^0\left(\widetilde V,\mathcal O_{\widetilde V}(D)\right)>0.$

From (\ref{eq4}) we need to show that
$K_{\widetilde V}+L\not\geq 0$ and $K_{\widetilde V}+M\not\geq 0.$
Suppose first that $K_{\widetilde V}+L\geq 0.$ Then
$$(K_{\widetilde V}+L)A_i=-1,\ \forall i\ \Longrightarrow\ K_{\widetilde V}+L-B\geq 0$$ and
$$
(K_{\widetilde V}+L-B)A_i'=-1,\ \forall i\ \Longrightarrow\ K_{\widetilde V}+L-B-C\geq 0.
$$
From $3K_{\widetilde V}+3L-3B-3C\geq 0$ and $3(L+M)\equiv 3(B+C)$ one gets
$3K_{\widetilde V}-3M\geq 0,$
i.e. $3K_{\widetilde V}-B-2C\geq 0.$
This implies the existence of an element in the linear system $|3K_V|$ having multiplicity
$> 1$ at each of the $12$ cusps of $V.$

Let $q_1,$ $q_2$ be the defining equations of $Q_1,$ $Q_2.$
The surface $V$ is given by equations $w^2=q_1,$ $t^2=q_2$ in the weighted projective
space $\m P(x^1,y^1,z^1,w^2,t^2).$
It is easy to see that the linear system of polynomials of degree $3$ has no element
with multiplicity $> 1$ at the cusps of $V$ (for instance using computer algebra).

The case $K_{\widetilde V}+M\geq 0$ is analogous.\Qed

\section{A surface with $K^2=9\chi$}

Based on the work of Prasad and Yeung \cite{PrYe}, Cartwright and Steger \cite{CaSt} constructed 
a fake projective plane $F$ with an automorphism $j$ of order $3$ such that $F/j$ is a  
surface of general type with $\chi=1,$ $p_g=0,$ $K^2=3$ and fundamental group $\m Z_{13}.$
This surface has a set of three $3$-divisible cusps (cf. \cite{Ke}).
Denote by $B$ the unit ball in $\m C\sp 2$ and let $P,$ $H$ and $G$ be the groups such that
$F=B/P,$ $F/j=B/H$ and the universal cover of $F/j$ is $B/G.$ 
We have the following commutative diagram, where the vertical arrows denote unramified $\m Z_{13}$
covers and the horizontal arrows denote $\m Z_3$ covers ramified over cusps.
$$
\begin{CD}\ B/{(G\cap P)}@>p>>B/G\\ @V 13:1 VV  @VV 13:1 V\\ \ B/P@>3:1 >> B/H
\end{CD}
$$

We show that the surface $S:=B/{(G\cap P)}$ is regular, hence $p_g(S)=p_g(B/G)$
and then the canonical map $\phi$ of $S$ factors through the triple cover $p$.
Since $G\cap P$ is the fundamental group of $S,$ the commutator quotient
$(G\cap P)/[G\cap P,G\cap P]$ is isomorphic to $H_1(S,\m Z)$. The first Betti number
$b_1(S)=2q(S)$ is the minimal number of generators of $H_1(S,\m Z)$ modulo elements of finite order.
Thus $q(S)=0$ if $(G\cap P)/[G\cap P,G\cap P]$ is finite.
This is shown in the Appendix, where we use computational GAP \cite{GAP4} data from
Cartwright and Steger to compute $G\cap P$.

Now notice that $\phi$ is not composed with a pencil. In fact otherwise the
canonical map of $B/G$ is composed with a pencil and then
$$39=K_{B/G}^2\geq 4\chi(\mathcal O_{B/G})-10=42,$$
from \cite[Theorem A]{Zu} (see also \cite[Corollary 3.4]{Ko}).

Finally we prove that if $d:=\deg(\phi)\ne 3,$
then $d=9$ and the canonical image $\phi(S)$ is a rational surface.
As in the proof of Proposition 4.1 of \cite{Be}, we have
$$9\chi(\mathcal O_{S})\geq K_{S}^2\geq d\ {\rm deg}(\phi(S))\geq nd(p_g(S)-2)$$
where $n=2$ if $\phi(S)$ is not ruled and $n=1$ otherwise.
This gives $d<6$ if $\phi(S)$ is not ruled and $d<12$ otherwise.
Since $d\equiv 0\ {\rm (mod}\ 3),$ then $d=6$ or $9$ and $\phi(S)$ is a rational surface.
If $d=6,$ the canonical map of $B/G$ is of degree $2$ and then $B/G$ has an involution.
But, as seen in the Appendix, the automorphism group of $B/G$ is the semidirect product
$\m Z_{13}:\m Z_3,$ so there is no involution on $B/G.$ 

\appendix
\section*{Appendix: GAP code}\label{appendix}

\begin{lstlisting}
#We are using data from
#http://www.maths.usyd.edu.au/u/donaldc/fakeprojectiveplanes/
#C18p3/C18p3-0-FP.gap
#namely the groups 'index9aFP', 'index3aFP' and
#the functions 'FundGp', 'AutGp'.

P:=index9aFP;
H:=index3aFP;
#B/P is a Fake projective plane.
#B/H is a quotient of the fake p.p. B/P by an order 3 automorphism.
#B/H is a surface with p_g=0 and K^2=3 having 3 cusps.

#Using the Cartwright-Steger function 'FundGp', we see that the
#fundamental group of B/H is "C13". We want to find the group G such
#that B/G is the universal cover of B/H. This is also computed by the
#function 'FundGp'.
#The following function is equal to 'FundGp' except that outputs G.

Grp:=function(G,FOList)
  local e1,e2,e3,GFO,G0,G0FCA;
  e1:=GeneratorsOfGroup(GammaBarFP);
  e2:=Concatenation(e1,List(e1,Inverse));
  Add(e2,One(G));
  e3:=ListX(e2,e2,\*);
  GFO:=Filtered(FOList,fo->fo in G);
  G0:=Group(ListX(e3,GFO,function(elt,fo) return fo^elt; end));
  Print(ForAll(GFO,fo->fo in G0),"\n");
  Print(IsNormal(G,G0),"\n");
  return G0;
end;

G:=Grp(H,FOList);
GP:=Intersection(G,P);

#The "commutative diagram":
StructureDescription(FactorGroup(P,GP)) = "C13";
StructureDescription(FactorGroup(H,G)) = "C13";
StructureDescription(FactorGroup(G,GP)) = "C3";
StructureDescription(FactorGroup(H,P)) = "C3";

#The first Betti number b_1(B/GP) = 0:
Index(GP,CommutatorSubgroup(GP,GP)) = 2916; #It is finite.

#The automorphism group of B/G:
StructureDescription(AutGp(G)) = "C13 : C3";
\end{lstlisting}

\bibliography{ReferencesRito}

\

\

\noindent Carlos Rito
\\ Departamento de Matem\' atica, Faculdade de Ci\^encias
\\ Rua do Campo Alegre 687, Apartado 1013
\\ 4169-007 Porto, Portugal
\\
{\it e-mail:} crito@fc.up.pt

\end{document}